\documentstyle{article}

\title{Solutions of Inverse Problems\\
for Variational Calculus}
\author{Mircea Neagu}
\date{}
\begin{document}
\maketitle
\begin{abstract}
In \S1 the author presents a short history of the problem studied in this
paper. \S2 introduces the notion of harmonic map between a Riemannian space
and a generalized Lagrange space, in a natural way. In \S3
is proved that for certain systems of differential or partial differential
equations, the solutions are harmonic maps, in the sense definite in \S1.
\S4 describes the main properties of the generalized Lagrange spaces constructed in \S3.
\end{abstract}
\section{Introduction}

\hspace{5mm}The problem of finding a geometrical structure of Riemannian type on a manifold
$M$ such that the orbits of an arbitrary vector field $X$ should be geodesics, was
analysed by Sasaki \cite{5}. Doing with this problem, Sasaki creates the well known
almost contact structures on a manifold of odd dimension, although the initial
problem rests open. After the introduction of
ge\-ne\-ra\-li\-zed Lagrange spaces by Miron \cite{2}, the same problem is resumed by
Udri\c ste \cite{6,7,8}. This succeded to discover a Lagrange structure on $M$,
dependent of the vector field $X$ and a (1,1)-tensor field, such that
the orbits of class $C^2$ should be geodesics.
Moreover, he formulates a more general problem \cite{8}, namely

1) There exist the structures of Lagrange type such that the solutions of
certain PDEs of order one should be {\it harmonic maps}?

2)What is a {\it harmonic map} between two spaces endowed with structures
of this type?\\
In this paper, the author tries to solve these problems. He will use the new
notion of harmonic map on a direction, offering thus a
partial answer to the Udri\c ste's questions.

\section{Harmonic maps on a direction}

\hspace{5mm} Let $(M^m,\;\varphi_{\alpha\beta}(a))$ be Riemannian manifold of
dimension $m$ and $(N^n,\;h_{ij}(x,y))$ generalized
Lagrange space of dimension $n$, where $a=(a^\alpha)_{\alpha=\overline{1,m}}$
are coordinates on $M$ and $(x,y)=(x^k,y^k)_{k=\overline{1,n}}$ are coordinates
on $TN$.

{\bf Remark.} On $M$, the coordinates are indexed by $\alpha,\beta,\gamma,
\ldots$ and, on $N$, respectively $TN$, the coordinates are indexed by
$i,j,k,\ldots$. Also, on $M\times N$, the first $m$ coordinates will be indexed
by $\alpha, \beta, \gamma\ldots$ and the last $n$ coordinates by $i,j,k,\ldots$.

Let $A\in{\cal X}(M)$ be an arbitrary vector field on $M$. If the manifold $M$
is connected, compact and orientable, we can define the $(\varphi,A,h)$-{\it
energy functional} or {\it the energy functional on the direction $A$} taking
$$
E^A_{\varphi h}:C^\infty (M,N)\to R,
$$
$$
E^A_{\varphi h}(f)=\displaystyle{{1\over 2}\int_M\varphi^{\alpha\beta}(a)
h_{ij}(f(a),f_*(A))f^i_\alpha f^j_\beta\sqrt\varphi da,}
$$
where $\displaystyle{f^i=x^i(f),\;f^i_\alpha ={\partial f^i\over\partial
a^\alpha},\;\varphi=det(\varphi_{\alpha\beta})}$ and $f_*:TM\to TN$ is the
differential of the map $f$.

{\bf Definition.}
The map $f\in C^\infty(M,N)$ is $(\varphi,A,h)${\it-harmonic} iff $f$ is a critical point
for the $(\varphi,A,h)$-energy functional $E^A_{\varphi h}$.

{\bf Remarks.} i) If $h_{ij}(x,y)=h_{ij}(x)$ is Riemannian metric, it recovers
the classical definition of a harmonic map between two Riemannian manifolds
\cite{1}.

ii) If we consider $M=[a,b]\subset R,\;\varphi_{11}=1$ and $A=d/dt$,
we obtain $C^\infty (M,N)=\{x:[a,b]\to N\vert\;x-C^\infty \hbox{differentiable}
\}\stackrel{\hbox{not}}=\Omega_{a,b}(N)$ and the $(1,d/dt,h)$-energy functional
should be
$$
E^{d/dt}_{1h}(x)=\displaystyle{{1\over 2}\int_a^bh_{ij}(x(t),\dot x(t))
{dx^i\over dt}{dx^j\over dt}dt,\;\forall x\in\Omega_{a,b}(N).}
$$
In conclusion, the $(1,d/dt,h)$-harmonic curves are exactly the geodesics of
the ge\-ne\-ra\-li\-zed Lagrange space $(N,h_{ij}(x,y))$ \cite {2}.

\section{The geometrical interpretation of the solutions of certain PDEs
of order one}

\hspace{5mm} Let $f\in C^\infty(M,N)$ be a smooth map and let be the global section
$$
\displaystyle{
\left.\delta f\stackrel{\hbox{not}}=f^i_\alpha da^\alpha\vert_a
\otimes{\partial\over\partial y^i}\right\vert_{f(a)}\in\Gamma(
T^*M\otimes f^{-1}(TN))}.
$$
On $M\times N$, let $T$ be one tensor of type $(1,1)$ with all
components equal to zero excepting $(T^i_{\alpha} )_{i=\overline{1,n}\\
\atop
\alpha=\overline{1,m}}$.
\medskip
Let be the system of partial differential equations
$$\displaystyle{
\delta f=T\;\hbox{expressed locally by}\;{\partial f^i\over\partial a^\alpha}
=T^i_\alpha(a,f(a)).
}\leqno(E)$$

If $(M,\varphi_{\alpha\beta})$ and $(N,\psi_{ij})$ are Riemannian
manifolds we can build a scalar pro\-duct on $\Gamma(T^*M\otimes f^{-1}(TN))$
putting
$<T,S>=\varphi^{\alpha\beta}(a)\psi_{ij}(f(a))T^i_\alpha S^j_\beta$, where
$\displaystyle{T=T^i_\alpha da^\alpha\otimes{\partial\over\partial y^i}}$ and
$\displaystyle{S=S^j_\beta da^\beta\otimes{\partial\over\partial y^j}}$.

In these conditions we can prove the following

{\bf Theorem.} {\it
If $(M,\varphi),(N,\psi)$ are Riemannian manifolds and $f\in C^\infty(M,N)$ is
solution of the system $(E)$, then $f$ is solution of the variational problem
asociated to the functional
${\cal L}_T:C^\infty(M,N)\backslash\{f\;\vert\;\exists\;a\in M\;\hbox{such that}\;
<\delta f,T>(a)=0\}\to R_+$,
$$
\displaystyle{{\cal L}_T(f)={1\over 2}\int_M{\Vert\delta f\Vert^2\Vert T\Vert^2
\over <\delta f,T>^2}\sqrt{\varphi}da={1\over 2}\int_M{\Vert T\Vert^2\over <
\delta f, T>^2}\varphi^{\alpha\beta}\psi_{ij}f^i_\alpha f^j_\beta\sqrt{\varphi}da}.
$$}

{\bf Proof.} In the space $\Gamma(T^*M\times f^{-1}(TN))$, the Cauchy
inequality for the scalar product $<,>$ holds.
It follows that we will have
$<T,S>^2\leq\Vert T\Vert^2\Vert S\Vert^2,\;\linebreak\forall\;
T,S\in\Gamma(T^*M\times\nolinebreak f^{-1}(TN)),$ with equality iff there exists
${\cal K}\in{\cal F}(M)$ such that $T={\cal K}S$. Consequently,
for every $f\in C^\infty (M,N)$, we will obtain
$$\displaystyle{{\cal L}(f)={1\over 2}\int_M{\Vert\delta f\Vert^2
\Vert T\Vert^2\over<\delta f,T>^2}\sqrt{\varphi}da\geq{1\over 2}\int_M\sqrt
{\varphi}da={1\over 2}Vol_\varphi(M)}.
$$
Now, if $f$ is solution of the system $(E)$, we will conclude
$\displaystyle{{\cal L}_T(f)={1\over 2}Vol_\varphi(M)}$. This means that
$f$ is a global minimum point for ${\cal L}_T$ $\Box$.

{\bf Remarks.} i) In certain particular cases of the system $(E)$, the
functional ${\cal L}_T$ becomes exactly a functional of type
$(\varphi,A,h)$-energy.

ii) The global minimum points of the functional ${\cal L}_T$ are solutions of
the system $\delta f={\cal K}T$, where ${\cal K}\in{\cal F}(M)$ not necessarily
with ${\cal K}=1$.

{\bf Examples.}

{\bf 1. Orbits}

For $M=([a,b],1)$ and $T=\xi\in\Gamma(x^{-1}(TN))$, the system
$(E)$ becomes
$$\displaystyle{{dx^i\over dt}=\xi^i(x(t)),\;x:[a,b]\to N},\leqno{(E_1)}$$
that is the system of orbits for $\xi$, and
the functional ${\cal L}_\xi$ is
$$\displaystyle{{\cal L}_\xi(x)={1\over 2}\int^b_a{\Vert\xi\Vert^2_\psi\over
<\xi,x_*(d/dt)>^2_\psi}\psi_{ij}{dx^i\over dt}{dx^j\over dt}dt}.
$$
Hence the functional ${\cal L}_\xi$ is a $(1,d/dt,h)$-energy
(see {\bf ii} of first remarks of this paper). The fundamental metric tensor
$h_{ij}:TN\backslash\{y\vert \xi^b(y)=0\}\to R$
is defined by
$$\displaystyle{
h_{ij}(x,y)={\Vert\xi\Vert^2_\psi\over<\xi,y>^2_\psi}\psi_{ij}(x)=\psi_{ij}(x)
\exp{\displaystyle{\left[
2\ln{\Vert\xi\Vert_\psi\over\vert<\xi,y>_\psi\vert}\right]} }
}.$$
This case is studied, in other way, by Udri\c ste in \cite{6}-\cite{8}.

{\bf 2. Pfaff systems}

For $N=(R,1)$ and $T=A\in\Lambda^1(T^*M)$, the system $(E)$
will become
$$df=A,\;f\in{\cal F}(M),\leqno{(E_2)}$$
that is a Pfaffian system, and the functional ${\cal L}_T$ reduces to
$$
\displaystyle{{\cal L}_A(f)={1\over 2}\int_M{\Vert A\Vert^2_\varphi\over
[f_*(A^\sharp)]^2}\varphi^{\alpha\beta}f_\alpha f_\beta\sqrt{\varphi}da},
$$
where
$\displaystyle{A^\sharp=\varphi^{\alpha\beta}A_\beta{\partial\over\partial a^\alpha}}.$
Hence the functional ${\cal L}_A$ is a
$(g,A^\sharp,h)$-energy, where
$$\displaystyle{
g_{\alpha\beta}(a)={1\over\Vert A\Vert^2_\varphi}\varphi_{\alpha\beta}(a)\;
\hbox{and}\;h:TR\backslash\{0\}\to R,\;
h(x,y)={1\over y^2}=e^{-2\ln\vert y\vert}.}$$

{\bf 3. Pseudolinear functions}

We suppose that $T^k_\beta (a,x)=\xi^{k}(x)A_\beta (a)$,
where $\xi^k$ is vector field on $N$ and $A_\beta$ is 1-form on $M$.
In this case the system $(E)$ will be
$$\displaystyle{{\partial f^k\over\partial a^\beta}=\xi^k(f)A_\beta (a)}
\leqno{(E_3)}
$$ and the functional ${\cal L}_T$ is expressed by
$$\displaystyle{{\cal L}_T(f)={1\over 2}\int_M{\Vert\xi\Vert^2_\psi\Vert A
\Vert^2_\varphi\over<\xi,f_*(A^\sharp)>^2}\varphi^{\alpha\beta}\psi_{ij}f^i_\alpha f^j_\beta
\sqrt\varphi da=
}$$
$$\displaystyle{={1\over 2}\int_Mg^{\alpha\beta}(a)h_{ij}(f(a),f_*(A^\sharp))
f^i_\alpha f^j_\beta\sqrt\varphi da},\;\hbox{where}
$$
$\displaystyle{g_{\alpha\beta}(a)={1\over\Vert A\Vert^2}
\varphi_{\alpha\beta}(a)\;\hbox{and}\;h_{ij}(x,y)=
{\Vert\xi\Vert^2\over<\xi,y>^2}\psi_{ij}(x)=
\exp{(-2\ln\vert<\xi,y>\vert)}\psi_{ij}(x).}$\\
It follows that the functional ${\cal L}_T$ becomes a $(g,A^\sharp,h)$-energy.

{\bf Remark.} Taking $M$ an open subset in $(R^n,\varphi=\delta)$ and
$N=(R,\psi=1)$, the system from the third example is
$$\displaystyle{{\partial f\over\partial a^\alpha}=\xi(a)A_\alpha(f(a))\;,\;
\forall\;\alpha=\overline{1,m}.}\leqno{(PL)}
$$
Supposing that $(grad\;f)(a)\ne 0\;,\;\forall\;a\in M$, the solutions of this
system are the well known {\it pseudolinear functions} \cite{4}. These functions
have the following property

$-$for every fixed point $x_0\in M$, the hypersurface of constant level
$$M_{f(x_0)}=\{x\in M\vert f(x)=f(x_0)\}$$ is totally geodesic
\cite {4} (i. e. the second fundamental form vanishes identically).

In conclusion, these pseudolinear functions are examples of
$(g,A^\sharp,h)$-harmonic functions on the Riemannian space
$\displaystyle{\left(M,g_{\alpha\beta}(a)={1\over\Vert A\Vert^2}
\delta_{\alpha\beta}\right)}$ with values into the ge\-ne\-ra\-lized Lagrange
space $\displaystyle{\left(R,h(x,y)={\xi(x)\over y}\right)}$. For example,
we have the following pseudolinear functions:\medskip

{\bf 3. 1.} $f(a)=e^{<v,a>+w}$, where $v\in M,\;w\in R$, is
solution for the system $(PL)$ with $\xi(a)=1$ and
$A(f(a))=f(a)v$.\medskip

{\bf 3. 2.} $\displaystyle{f(a)={<v,a>+w\over<v',a'>+w'}}$, where $v,v'\in M
\;;\;w,w'\in R$, is solution for $(PL)$ with $\displaystyle{\xi(a)={1\over<v',a>+w'}}$
and $A(f(a))=v-f(a)w$.

{\bf Remark.} A geometrical interpretation of the solutions of the system $(E)$,
in the general case, when the tensor $T$ is expressed by
$T^i_\alpha(a,x)=\sum^t_{r=1}\xi^i_r(x)A^r_\alpha(a)$, where $\{\xi_r\}_{r=
\overline{1,t}}\subset{\cal X}(N)$ is a family of vector fields on $N$ and
$\{A^r\}_{r=\overline{1,t}}\subset\Lambda^1(T^*M)$ is a family of covector
fields on $M$, will be treated, in other sense, in a subsequent paper.

\section{The geometry asociated to PDEs of order one}

\hspace{5mm}For beginning, we remark that, in all above cases, the solutions
of class $C^2$ of the system $\delta f=T$
becomes $(\varphi,A,h)$-harmonic maps, in the sense definite in this paper.
Moreover, the generalized Lagrange structures constructed above are of type
$(M^n,e^{2{\sigma(x,y)}}\gamma_{ij}(x))$, where $\sigma :TM\backslash\{
\hbox{Hyperplane}\}\to R$ is a smooth function.
Now, we assume that a generalized Lagrange space $(M^n,g_{ij}(x,y))$
satisfies the following axioms:

a. 1. The fundamental tensor field $g_{ij}(x,y)$ is of the form
$$
g_{ij}(x,y)=e^{2\sigma(x,y)}\gamma_{ij}(x).
$$

a. 2. The space is endowed with the non-linear connection
$$
N^i_j(x,y)=\Gamma^i_{jk}(x)y^k,
$$
where $\Gamma^i_{jk}(x)$ are the Christoffel symbols for the Riemannian
metric $\gamma_{ij}(x)$.

Under these assumptions, our space verifies a constructive-axiomatic formulation of General Relativity
due to Ehlers, Pirani and Schild \cite{2}. This space represents a convenient
relativistic model, since it has the same conformal and projective properties as the
Riemannian space $(M,\gamma_{ij})$.

Developing the formalism presented in \cite{2}, \cite{3} and denoting the curvature tensor
field of the metric $\gamma_{ij}(x)$ by $r^i_{jkl}$, the following Maxwell's
equations hold
$$\left\{\begin{array}{lll}
\displaystyle{F_{ij\vert k}+F_{jk\vert i}+F_{ki\vert j}=-\sum_{(ijk)}g_{ip}
r^h_{qjk}{\partial\sigma\over\partial y^h}y^py^q},\\
F_{ij}\vert_k+F_{jk}\vert_i+F_{ki}\vert_j=-(f_{ij\vert k}+
f_{jk\vert i}+f_{ki\vert j}),\\
f_{ij}\vert_k+f_{jk}\vert_i+f_{ki}\vert_j=0,
\end{array}\right.
$$
where the electromagnetic tensors $F_{ij}$ and $f_{ij}$ are
$$\displaystyle{
F_{ij}=\left(g_{ip}{\delta\sigma\over\delta x^j}-g_{jp}{\delta\sigma\over
\delta x^i}\right)y^p\;,\;f_{ij}=\left(g_{ip}{\partial\sigma\over\partial y^j}-
g_{jp}{\partial\sigma\over\partial y^i}\right)y^p
}.$$

We will use the following notations:
$r_{ij}=r^k_{ijk}$, $r=\gamma^{ij}r_{ij}$,
$\displaystyle{
{\delta\over\delta x^i}={\partial\over \partial x^i}-
N^j_i{\partial\over\partial y^j}}$, $\sigma^{\scriptscriptstyle {H}}=
\displaystyle{\gamma^{kl}{\delta\sigma\over\delta x^k}{\delta\sigma\over\delta
x^l},\;\sigma^{\scriptscriptstyle {V}}=\gamma^{ab}{\partial\sigma\over\partial
y^a}{\partial\sigma\over\partial y^b},\;\overline\sigma=\gamma^{ij}\sigma_{ij},
\dot\sigma=\gamma^{ab}\dot\sigma_{ab}}$,
$$\hbox{where}\quad\left\{
\begin{array}{ll}\medskip
\displaystyle{
\sigma_{ij}={\delta\sigma\over\delta x^i}\vert_j+{\delta\sigma\over\delta x^i}{\delta\sigma\over
\delta x^j}-{1\over 2}\gamma_{ij}\sigma^{\scriptscriptstyle H}}\\
\displaystyle{
\dot\sigma_{ab}=\left.{\partial\sigma\over\partial y^a}\right\vert_b+
{\partial\sigma\over\partial y^a}{\partial\sigma\over\partial y^b}-{1\over 2}
\gamma_{ab}\sigma^{\scriptscriptstyle V}}.\\
\end{array}\right.
$$
The Einstein's equations of the space $(M,g_{ij}(x,y))$ take the form
$$\left\{
\begin{array}{ll}\medskip
\displaystyle{r_{ij}-{1\over 2}r\gamma_{ij}+t_{ij}={\cal K}
T^{\scriptscriptstyle H}_{ij}}\\
(2-n)(\dot\sigma_{ab}-\dot\sigma\gamma_{ab})={\cal K}T^{\scriptscriptstyle V}
_{ab},\\
\end{array}\right.
$$
where $T^{\scriptscriptstyle H}_{ij}$ and $T^{\scriptscriptstyle V}_{ab}$ are
the $h$- and the $v$-components of the energy momentum tensor field, ${\cal K}$
is the gravific constant and
$$
t_{ij}=(n-2)(\gamma_{ij}\overline\sigma-\sigma_{ij})+\gamma_{ij}r_{st}y^s
\gamma^{tp}{\partial\sigma\over\partial y^p}+{\partial\sigma\over\partial y^i}
r^a_{tja}y^t-\gamma_{is}\gamma^{ap}{\partial\sigma\over\partial y^p}
r^s_{tja}y^t.
$$

Consequently, in certain particular cases,
it is posible to build a generalized Lagrange geometry (in Miron's
sense) naturally attached to a system of partial di\-ffe\-ren\-ti\-al equations.
This idea was suggested by Udri\c ste in private discussions.

{\bf Open problem.} Because the generalized Lagrange structure constructed
in this paper is not unique, it arises a natural question:

-Is it possible to build a unique generalized Lagrange geometry naturally
asociated to a given PDEs system?

\begin{tabbing}
\kill\\
\end{tabbing}
\begin{center}
University POLITEHNICA of Bucharest\\
Department of Mathematics I\\
Splaiul Independentei 313\\
77206 Bucharest, Romania\\
e-mail:mircea@mathem.pub.ro
\end{center}

\end{document}